\documentclass[a4paper,11pt]{article}

\usepackage{latexsym}
\usepackage{amsmath,amsfonts,amssymb}
\usepackage{graphicx,color,epsfig}
\usepackage{graphics,epsf,psfrag}

\def\C{\mathbb{C}}
\def\N{\mathbb{N}}

\def\R{\mathbb{R}}
\def\T{\mathbb{T}}
\def\Z{\mathbb{Z}}

\def\leq{\leqslant}
\def\geq{\geqslant}
\def\eps{\epsilon}

\numberwithin{equation}{section}

\newtheorem{Thm}{Theorem}[section]
\newtheorem{Lem}[Thm]{Lemma}
\newtheorem{Pro}[Thm]{Proposition}
\newtheorem{Cor}[Thm]{Corollary}

\begin{document}

\title{Landau damping in the Kuramoto model}
\author{Bastien Fernandez$^{1,2}$, David G\'erard-Varet$^3$ and Giambattista Giacomin$^1$}
\maketitle

\footnotetext[1]{Universit\'e Paris Diderot, Sorbonne Paris Cit\'e,  Laboratoire de Probabilit{\'e}s et Mod\`eles Al\'eatoires, UMR 7599 Universit\'e Pierre et Marie Curie, CNRS, F- 75205 Paris, France}

\footnotetext[2]{On leave from Centre de Physique Th\'{e}orique, CNRS - Aix-Marseille Universit\'e - Universit\'e de Toulon, Campus de Luminy, 13288 Marseille CEDEX 09 France}

\footnotetext[3]{Universit\'e Paris Diderot, Sorbonne Paris Cit\'e, Institut de Math\'ematiques de Jussieu -- Paris Rive Gauche, UMR 7586,
 F- 75205 Paris, France}

\begin{abstract}
We consider the Kuramoto model of globally coupled phase oscillators in its continuum limit, with individual frequencies drawn from a distribution with density of class $C^n$ ($n\geq 4$). A criterion for linear stability of the uniform stationary state is established which, for basic examples of frequency distributions, is equivalent to the standard condition on the coupling strength in the literature. We prove that, under this criterion, the  Kuramoto order parameter, when evolved under the full nonlinear dynamics, asymptotically vanishes (with polynomial rate $n$) for every trajectory issued from sufficiently small $C^n$ perturbation. The proof uses techniques from the Analysis of PDEs and closely follows recent proofs of the nonlinear Landau damping in the Vlasov equation and Vlasov-HMF model.
\end{abstract}

{\small
\noindent
2010 \textit{Mathematics Subject Classification:} 35Q84, 82C44, 45D05, 35Q92,  92B25

\noindent
\textit{Keywords:} Kuramoto model, Disordered mean field rotators, Landau damping, Volterra integral equation
}

\section{Introduction}
The Kuramoto model is the archetype of nonlinear dissipative systems of coupled oscillators. Introduced about forty years ago to mimic simple chemical instabilities \cite{K75,K84}, this model has since been applied to a large palette of systems in various domains, from Physics to Biology, to Ecology and to Social Sciences, see \cite{ABP-VRS05,S00} for some examples of application. 

In its basic form, this model considers a population of $N$ oscillators ($N\in\N$, supposedly large) that are characterized by their phase on the unit circle, {\sl viz.}\ $\theta_i\in\T^1=\R/2\pi\Z$, and whose dynamics is governed by the following system of globally coupled ODEs
\[
\frac{d\theta_i}{dt}=\omega_i+\frac{K}{N}\sum_{j=1}^N\sin (\theta_j-\theta_i),\ \forall i\in\{1,\cdots,N\},t>0.
\]
The individual, time-independent, frequencies $\omega_i\in\R$ are randomly drawn using a probability density $g$ (which was assumed to be symmetric and unimodal in the original formulation). The parameter $K\in\R^+$ measures the interaction strength. 

Many variations on the theme can now be found in the voluminous literature on this subject, including singular or bimodal densities $g$ \cite{CCHKJ14,CHJK12,HHK10,MBSOSA09}, interactions involving several Fourier modes \cite{cf:daido} and random perturbations of the deterministic dynamics by additive noise \cite{ABP-VRS05,cf:GLP,SM91}, see also \cite{IPM-CS13,MP11} for recent generalizations. Here, we shall be specially interested in the continuum formulation of the original model that is obtained at limit $N\to\infty$. 

The main features of the Kuramoto phenomenology were identified in the early studies, using both numerical simulations and the analysis of stationary states \cite{K75,K84}. In few words, the dynamics is as follows (we refer to \cite{S00} for a concise yet more detailed account). While the population behavior remains fully incoherent for sufficiently small interaction strengths, a regime of partial synchronization, in which all oscillators with sufficiently small frequency are locked together, takes place when $K$ increases beyond some threshold $K_c$. (Increasing $K$ further, full synchronization can eventually be achieved, provided that $g$ has compact support.)

The transition at $K_c$ can be identified using the order parameter $\text{R}_N(t)$ associated with each trajectory $t\mapsto\{\theta_i(t)\}_{i=1}^N$ and defined by 
\[
\text{R}_N(t)=\sum_{j=1}^Ne^{i\theta_j(t)}.
\]
In fact, no matter what the initial condition is, the quantity $\text{R}_N(t)$ asymptotically vanishes when $t\to +\infty$, as long as $K\in [0,K_c)$. For $K>K_c$, the quantity $|\text{R}_N(t)|$ tends to a positive value, and this limit increases with $K$ (and reaches 1 when full synchronization occurs). These findings were confirmed later on by thorough investigations on the linear stability of stationary states, both for the finite dimensional system and its continuum limit \cite{MS05,MS07,SMM92}.  

For symmetric and unimodal densities $g$, the threshold is given by $K_c=\frac2{\pi g(0)}$ and the transition corresponds to a supercritical pitchfork bifurcation of the order parameter \cite{ABP-VRS05,K75,K84,S00}. In other cases, the expression of $K_c$ and/or the bifurcation type may differ. For instance, the bifurcation at $\frac2{\pi g(0)}$ is subcritical when $g''(0)>0$ and see \cite{MBSOSA09} for the expression of $K_c$ in the case of symmetric bi-Cauchy distribution densities. 

In spite of accurate results for the linearized dynamics, to the best of our knowledge, full proofs of this phenomenology ({\sl i.e.}\ with nonlinear terms  included) remained to be provided. Of note, assuming Gaussian or Cauchy distribution, a control of the nonlinear dynamics via center manifold reduction is exposed in \cite{C13,CN11}, and a complete bifurcation analysis in the four dimensional invariant subspace based on the so-called Ott-Antonsen ansatz \cite{OA09} is reported in \cite{MBSOSA09} (for symmetric bi-Cauchy distributions as previously mentioned). Besides, when all individual frequencies are equal (or, equivalently, when the frequency distribution is supported at the origin $\omega=0$), complete asymptotic synchrony has been proved to hold for every $K>0$ \cite{HHK10}.

In this paper, we treat from a rigorous mathematical viewpoint, the complete nonlinear dynamics in the incoherent regime of the continuum limit of the Kuramoto model. By identifying a stability condition that holds for every sufficiently regular density $g$, we prove that for every sufficiently regular and small initial perturbation
of the uniform stationary state, the order parameter of the associated trajectory vanishes as $t\to +\infty$. We also show that the same conclusion holds for any sufficiently regular perturbation (not necessarily small), provided that the interaction strength is small enough. To some extent, these results confirm that the phenomenology mentioned above holds for a large class of frequency densities and initial perturbations of the infinite-dimensional system.

Instead of using Dynamical System techniques (which were used in previous mathematical approaches), our analysis relies on PDEs methods and more precisely, on the approach to Landau damping in the Vlasov-HMF model in \cite{FR14}. In particular, our main result is a uniform in time control of the solution regularity, which implies (polynomial) decay in Fourier space (the more regularity, the faster the decay), from where the property $\text{R}(t)\to 0$ follows suit. 

The analogy between the Kuramoto model and the Vlasov equation -- especially that the Kuramoto order parameter $\text{R}(t)$ is the analogue of the so-called force (hence the term Landau damping for the property $\text{R}(t)\to 0$) -- had already been employed in the literature, in particular, for the linear stability analysis, see {\sl e.g.}\ \cite{SMM92}. In both cases, the linearized dynamics can be regarded as a Volterra equation (of the second kind) that can be controlled under appropriate conditions on the associated convolution kernel (see \cite{V10} for a pedagogical exposition of the Vlasov equation analysis). Our stability criterion is actually the analogue of the stability criterion for the Vlasov equation. Furthermore, the analogy had also been used to prove that the dynamics of the continuum limit is a suitable approximation of the finite dimensional system dynamics on every finite-time interval \cite{L05}. 

A complete approach to nonlinear Landau damping in the Vlasov(-Poisson) equation was first established in \cite{MV11}, and recently improved in \cite{BMM13}. In particular, one of the prowesses in \cite{MV11} was to deal with singular potentials. Instead, the potential in Vlasov-HMF only consists of the first Fourier mode, as in the Kuramoto model; hence the bootstrap argument in \cite{FR14} suffices for our purpose. In analogy with \cite{FR14},  our results can be extended to variations of the Kuramoto model where the interaction consists of a finite number of harmonics \cite{cf:daido}. 

To our knowledge, this paper is the first application of nonlinear damping methods to the global phenomenology of coupled oscillator systems with dissipative dynamics. We hope that this transfer of techniques can be extended to other models, beyond the Kuramoto setting.

The paper is organized as follows. The PDE under consideration and associated quantities are introduced in Section 2. Section 3 contains the main results of the paper and related comments, especially on the relationships between our stability condition and previous ones in the literature. In Section 4, we obtain a Volterra equation for the (rescaled) order parameter and we provide a control of the solution, depending on the input term. Section 5 completes the proof of Landau damping by reporting the analysis of the full nonlinear equation based on a bootstrap argument. 
 
\section{Continuum limit of the Kuramoto model}
When restricted to absolutely continuous distributions of the cylinder $\T^1\times\R$, the continuum limit of the Kuramoto model is given by the following PDE \cite{L05,SM91} for densities $\rho(\theta,\omega)\geq 0$ such that $\int_{\T^1}\rho(\theta,\omega)d\theta=1$ for all $\omega\in\R$,
\begin{equation}
\partial_t\rho(\theta,\omega)+\omega\partial_\theta\rho(\theta,\omega) +\partial_\theta\left(\rho(\theta,\omega)V(\theta,\rho)\right)=0,\ \forall (\theta,\omega)\in \T^1\times\R,\ t>0,
\label{BASICDYNAM}
\end{equation}
where the potential $V$ is defined by 
\[
V(\theta,\rho)=K\int_{\T^1\times\R}\sin (\vartheta-\theta)\rho(\vartheta,\omega)g(\omega)d\vartheta d\omega.
\]
Up to a rescaling of time and frequencies in equation \eqref{BASICDYNAM}, the parameter $K$ could be absorbed in the probability density $g$ (which is arbitrary at this stage). Hence, similarly to the Vlasov equation, this density is the only relevant parameter in the dynamics. Nevertheless, we keep the current formulation of equation \eqref{BASICDYNAM}, in agreement with the original Kuramoto conjecture on the dynamics dependence on the interaction strength for a given probability density. 

Throughout the paper, we shall need the following notations. Given a real function $u=u(\omega)$ of class $C^n$, consider its Sobolev norm $\|\cdot\|_{{\cal H}^n}$ defined by 
\[
\|u\|_{{\cal H}^n}^2=\sum_{k=0}^n\|\langle \omega\rangle u^{(k)}\|_{L^2(\R)}^2,
\]
where $\langle \omega\rangle=\sqrt{1+\omega^2}$ for all $\omega\in\R$ and $u^{(k)}$ denotes the $k$th derivative. For any $L^1$ real function $u$, let 
\[
\widehat{u}(\tau)=\int_\R u(\omega)e^{-i\tau\omega}d\omega,\ \forall \tau\in\R,
\]
be its Fourier transform.

Assuming that the initial density $\rho(0,\theta,\omega)$ at time 0 is of class $C^n$ ($n\in\N$) on the cylinder, equation \eqref{BASICDYNAM} has a unique global solution $t\mapsto \rho(t,\theta,\omega)$ defined on $\R^+$ which is of class $C^n$ and $\rho(t,\theta,\omega)$ remains a normalised density at all times (global existence can be established by standard arguments, for instance, by applying the method of characteristics \cite{L05}).  

The uniform density $\rho(\theta,\omega)=\frac{1}{2\pi}$ is an obvious stationary solution of equation \eqref{BASICDYNAM}. In order to study the dynamics of perturbations to this uniform steady state, we consider the following decomposition
\begin{equation}
\rho(t,\theta,\omega)=\frac{1}{2\pi}+\text{r}(t,\theta,\omega),\ \forall t\in\R^+.
\label{DECOMPER}
\end{equation}
(Notice that we must have $\text{r}(t,\theta,\omega)\geq -\frac{1}{2\pi}$ and $\int_{\T^1}\text{r}(t,\theta,\omega)d\theta=0$ for all $\omega\in\R$). Equation \eqref{BASICDYNAM} implies that $\text{r}$ must be the unique global $C^n$ solution - which exists for any $C^n$ initial condition $r(0)$ - of the  following equation
\begin{equation}
\label{eq:r}
\partial_t\text{r}(\theta,\omega)+\omega\partial_\theta \text{r}(\theta,\omega) +\partial_\theta\left(\left(\frac{1}{2\pi}+ \text{r}(\theta,\omega)\right)V(\theta,\text{r})\right)=0,\ \forall (\theta,\omega)\in \T^1\times\R,\ t>0.
\end{equation}

Similarly to as above for real functions, we shall employ the Sobolev norm of a function $u$ of class $C^n$ on the cylinder, namely the quantity (which we denote using the same symbol as before)
\[
\|u\|_{{\cal H}^n}^2=\sum_{k_\theta,k_\omega\geq 0,\ k_\theta+k_\omega\leq n}\|\langle \omega\rangle\partial_\theta^{k_\theta}\partial_\omega^{k_\omega}u\|_{L^2(\T^1\times\R)}^2.
\]

\section{Main results}
\subsection{Landau damping in the Kuramoto model}
In this section, we formulate and comment our result on the asymptotic behavior as $t\to +\infty$ of the order parameter $\text{R}(t)$ defined by 
\[
\text{R}(t)= \int_{\T^1\times\R} \text{r}(t,\theta,\omega)g(\omega)e^{-i\theta}d\theta d\omega ,
\] 
which obviously satisfies $|\text{R}(t)|\leq 1$. The behavior of $\text{R}(t)$ is given in the following statement. Let 
\[
\Pi^-=\{x+iy:\, x\in\R, \, y\in\R^-\},
\]
be the lower half plane of complex numbers with negative imaginary part and notice that for any function $G\in L^1(\R)$, the integral $\int_{\R^+}G(t)e^{-i\omega t}dt$ is finite for every $\omega\in \Pi^-$.
\begin{Thm}
Assume that $g$ is of class $C^n(\R)$ for some $n\geq 4$ and satisfies the following conditions 
\begin{equation}
\|g\|_{{\cal H}^n}<+\infty,\quad \widehat{g}\in L^1(\R^+)\quad\text{and}\quad\int_{\R^+}t^n|\widehat{g}(t)|dt<+\infty. 
\label{CONDG}
\end{equation}
Then, for every $K\geq 0$ such that 
\begin{equation}
1-\frac{K}{2}\int_{\R^+}\widehat{g}(t)e^{-i\omega t}dt\neq 0,\ \forall \omega\in \Pi^-,
\label{STABCRIT}
\end{equation}
there exists $\epsilon_K>0$ such that for any initial probability density $\frac1{2\pi}+\text{\rm r}$ where $\text{\rm r}$ is of class $C^n(\T^1\times\R)$ and such that $\|\text{\rm r}\cdot g\|_{{\cal H}^n}\leq \epsilon_K$, the order parameter associated with the solution of equation \eqref{BASICDYNAM} has the following asymptotic behavior
\[
\text{\rm R}(t)=O(t^{-n}).
\]
\label{MAINRES}
\end{Thm}

For the proof, see beginning of section \ref{S-PROOF} below. Theorem \ref{MAINRES} calls for the following series of comments. 
\begin{itemize}
\item The conclusion $\text{\rm R}(t)=O(t^{-n})$ is only non-trivial for initial perturbations such that $\text{R}(0)=\int_{\T^1\times\R} \text{r}(0,\theta,\omega)g(\omega)e^{-i\theta}d\theta d\omega\neq 0$. Otherwise, we would have $\text{R}(t)=0$ and $\text{r}(t,\theta,\omega)=\text{r}(0,\theta-t\omega,\omega)$ for all $t\in\R^+$ (consequence of the uniqueness 
for the Volterra equation \eqref{FULLEQ} below with $\widehat{p}_1(0,t)=0$ for all $t\in\R^+$). 
\item The assumption $\|g\|_{{\cal H}^n}<+\infty$ implies that the condition $\|\text{r}\cdot g\|_{{\cal H}^n}\leq \epsilon_K$ holds in particular for any perturbation $\text{r}$ such that the norm
\[
\max_{k_\theta,k_\omega\geq 0,\ k_\theta+k_\omega\leq n}\|\partial_\theta^{k_\theta}\partial_\omega^{k_\omega}\text{r}\|_{L^\infty(\T^1\times\R)},
\]
is sufficiently small.
\item Theorem \ref{MAINRES} also reveals that the decay $O(t^{-n})$ of the order parameter results from some interplay between the ${\cal H}^n$ regularity of $g$ and $\text{r}$; higher regularity implies stronger decay. This rate is likely to be optimal since the one of the linearized dynamics is already given by the regularity of $\text{r}$ \cite{SMM92} (see also Proposition \ref{GENCONTROL}). On the other hand, the rate does not depend on $K$. 
\item The minimal required regularity $n\geq 4$ is a by-product of our technical estimates. We do not know if Landau damping holds for every, say $C^1$ perturbations, or if there can be arbitrarily small unstable (not $C^4$) perturbations for which $|\text{R}(t)|\geq \delta>0$ uniformly in time.  
\end{itemize}
In addition, the $K$-dependent smallness condition $\|\text{r}\cdot g\|_{{\cal H}^n}\leq \epsilon_K$ on the size of perturbations can be justified by the fact that this size should vanish when $K$ approaches the value where the criterion \eqref{STABCRIT} fails  (especially in the case of a subcritical bifurcation where incoherent and partially locked stationary states co-exist for $K<K_c$). However, a closer look at proof of Theorem \ref{MAINRES} shows that any perturbation in ${\cal H}^n$ can be made admissible provided that $K$ is sufficiently small, as now claimed.
\begin{Pro}
Under the same conditions on $g$ as in Theorem \ref{MAINRES}, for any initial probability density $\frac1{2\pi}+\text{\rm r}$ with $\|\text{\rm r}\cdot g\|_ {{\cal H}^{n}}<+\infty$, there exists $K_r>0$ such that the conclusion of Theorem \ref{MAINRES} holds for every $K\in [0,K_r)$.
\label{LARGEPERT}
\end{Pro}
The proof of this statement is given in section \ref{S-LARGEPERT}.

Finally, and as for the Vlasov equation \cite{FR14,V10}, a direct consequence of Landau damping is the weak convergence of the solution $\text{r}(t)$ of \eqref{eq:r}
to a solution of the free transport equation 
\[
\partial_t\text{r}(\theta,\omega)+\omega\, \partial_\theta \text{r}(\theta,\omega)=0.
\]
To see this, given $t\in\R^+$ and a function $u$ of the cylinder, consider the Galilean change of variables $T^t$ defined by
\[
T^tu(\theta,\omega)=u(\theta+t\omega,\omega),\ \forall (\theta,\omega)\in \T^1\times\R.
\]
\begin{Cor}
Under the conditions of Theorem \ref{MAINRES} (resp.\ Proposition \ref{LARGEPERT}), for any initial probability density $\frac1{2\pi}+\text{\rm r}$ with $\|\text{\rm r}\cdot g\|_ {{\cal H}^{n}}<\epsilon_K$, (resp.\ $\|\text{\rm r}\cdot g\|_ {{\cal H}^{n}}<+\infty$), there exists a function $\text{\rm r}_\infty$ of the cylinder with $\|\text{\rm r}_\infty \cdot g\|_{{\cal H}^{n-2}}<+\infty$ such that the perturbation $\text{\rm r}(t)$, associated with the solution of equation \eqref{BASICDYNAM}, has the following asymptotic behavior
\[
\lim_{t\to+\infty}\|(T^t\text{\rm r}(t)-\text{\rm r}_\infty)\cdot g\|_{{\cal H}^{n-2}}=0.
\]
\label{WEAKCONV}
\end{Cor}
We provide a proof in section \ref{S-WEAKCONV}.

\subsection{Analysis of the stability criterion}\label{S-ANALSTABCOND}
This section reports some comments on the conditions \eqref{CONDG} and \eqref{STABCRIT} in Theorem \ref{MAINRES} and their relationships with previous stability conditions in the literature.

The integrability conditions on $\widehat{g}$ in \eqref{CONDG} can be granted by imposing enough regularity on $g$. In particular, since the square of $\int _{\R}\vert g^{(k)}\vert  $ is bounded by
$\int_{\R}\langle \omega \rangle^{-2} \int_{\R} \langle \omega \rangle^{2} \vert g^{(k)}\vert^2$ we have
\begin{equation}
t^k|\widehat{g}(t)|\leq \sqrt{\pi}\|g\|_{{\cal H}^k},\ \forall t\in\R^+,k\in\N,
\label{CONTFOUR}
\end{equation}
and therefore the desired integrability holds provided that $g$ is of class $C^{n+2}$ and $\|g\|_{{\cal H}^{n+2}}<+\infty$. 

\noindent
The conditions \eqref{CONDG} hold for all $n\geq 4$ for the (density of the) Cauchy distribution 
\[
g_\Delta(\omega)=\frac{\Delta}{\pi(\omega^2+\Delta^2)},
\]
and therefore, for any finite convex combination $\sum_i \alpha_i g_{\Delta_i}(\cdot-\omega_i)$. The same property holds for the Gaussian distribution $g_{\sigma}(\omega)=\frac1{\sigma\sqrt{2\pi}}e^{-\frac12\left(\frac{\omega}{\sigma}\right)^2}$. These examples, especially the densities $g_\Delta$ and $\frac12(g_\Delta(\cdot+\omega_0)+g_\Delta(\cdot-\omega_0))$, have been extensively considered in the literature, see {\sl e.g.}\ \cite{ABP-VRS05,IPM-CS13,K75,MBSOSA09,MP11,SM91}.
\medskip

Condition \eqref{STABCRIT} is the analogue of the stability criterion of the Vlasov equation \cite{V10}. It appears to be optimal, as least, as far as linear stability is concerned. To see this, assume the existence of $\omega_0$ with $\text{Im}(\omega_0)<0$ such that 
\[
\frac{K}{2}\int_{\R^+}\widehat{g}(t)e^{-i\omega_0 t}dt=1.
\]
Then a direct calculation based on Fubini's Theorem shows that we have 
\[
\frac{K}{2}\int_\R\frac{g(\omega)}{i(\omega_0+\omega)}d\omega=1,
\]
and this condition implies that the linearized Kuramoto equation 
\[
\partial_t\text{\rm r}(\theta,\omega)+\omega\partial_\theta \text{\rm r}(\theta,\omega) +\frac{1}{2\pi}\partial_\theta V(\theta,\text{\rm r})=0,
\]
has a solution (see section 3 in \cite{SM91})
\[
{\text r}(t,\theta,\omega)=\frac{Ae^{i\theta}}{i(\omega+\omega_0)}e^{i\omega_0t},\ \forall t\in \R^+,
\]
(with $A\neq 0$), whose order parameter $\text{R}(t)=\frac{4\pi A}{K}e^{i\omega_0t}$ grows exponentially with $t$. This exponential instability can be alternatively exhibited in the integral formulation of the dynamics, see Lemma \ref{OPTIMAL} below (for the case of a Volterra equation with arbitrary kernel). In addition, if \eqref{STABCRIT} fails for some $\omega_0$ in the real axis, one can prove that the order parameter $\text{R}$ associated to the solution of linearized equation cannot belong to $L^1(\R^+)$; hence the term 'optimal' for this criterion.

Furthermore, as for the Vlasov equation, one can obtain more explicit stability conditions than \eqref{STABCRIT}. To that goal, we first need the following considerations. Given a complex valued function $G\in L^1(\R^+)$, let 
\begin{equation}
D_G(\omega)=1-\int_{\R^+}G(t)e^{-i\omega t}dt,\ \forall \omega\in\R.
\label{DEFDG}
\end{equation}
The function $D_G$ is continuous on $\R$ and the Riemann-Lebesgue lemma implies $\lim\limits_{\omega\to \pm\infty}D_G(\omega)=1$. Extending $D_G$ by continuity to $\overline{\R}$, the expression
\[
\gamma_G=\{D_G(\omega)\}_{\omega\in\overline{\R}},
\]
defines a closed path in the complex plane. Assuming $D_G|_{\R}\neq 0$, let $\text{Ind}_G(0)\in\Z^+$ be the index (winding number) of $0$ with respect to $\gamma_G$. We have the following statement whose proof is given in section \ref{S-EQUISTAB} below.

\begin{Lem}
Assume that $G\in L^1(\R^+)$, $\int_{\R^+}t|G(t)|dt<+\infty$ and $D_G|_{\R}\neq 0$. Then we have $D_G|_{\Pi^-}\neq 0$ iff  $\text{\rm Ind}_G(0)=0$.
\label{EQUISTAB}
\end{Lem}
Under the condition \eqref{CONDG}, the function $\widehat{g}$ satisfies the first two assumptions of this Lemma and we have $\lim\limits_{\omega\to \pm\infty}D_{\frac{K}2\widehat{g}}(\omega)=1$. Therefore, to make sure that condition \eqref{STABCRIT} holds for a density $g$ satisfying \eqref{CONDG}, it suffices to ensure that the real part of $D_{\frac{K}2\widehat{g}}(\omega)$ is positive at every $\omega\in\R$ for which the imaginary part of $D_{\frac{K}2\widehat{g}}(\omega)$ vanishes. Using the expression 
\[
\int_{\R^+}\widehat{g}(t)e^{-i\omega t}dt=
\pi g(\omega)+ i\int_{\R^+}\frac{g(\omega-\sigma)-g(\omega+\sigma)}{\sigma}d\sigma,\forall \omega\in\R,
\]
we obtain the following sufficient condition for stability, which already appeared in \cite{IPM-CS13,MP11} and which is the analogue of the Penrose criterion \cite{V10} for the Vlasov equation
\begin{equation}
\int_{\R^+}\frac{g(\omega-\sigma)-g(\omega+\sigma)}{\sigma}d\sigma=0\quad\Longrightarrow\quad K<\frac{2}{\pi g(\omega)}.
\label{PENROSE}
\end{equation}
When $g$ is symmetric and unimodal, the integral only vanishes for $\omega=0$ and,  as observed in \cite{IPM-CS13}, this criterion reduces to the original inequality  $K<\frac{2}{\pi g(0)}$. Moreover, criterions \eqref{STABCRIT} and \eqref{PENROSE} are equivalent and optimal for the full nonlinear dynamics, in the sense that partially locked stationary solutions with $|\text{R}|>0$ are well-known to exist for $K>\frac{2}{\pi g(0)}$ in this case. 

Furthermore, criterions \eqref{STABCRIT} and \eqref{PENROSE} are also equivalent and optimal for the bi-Cauchy density $g_{\Delta, \omega_0}(\cdot):=\frac12(g_\Delta(\cdot+\omega_0)+g_\Delta(\cdot-\omega_0))$. Indeed, explicit calculations yield
\[
D_{\frac K2 \widehat g_{\Delta, \omega_0}} (\omega) \, =\,
1 -\frac K2 \left(\frac{\Delta + i \omega }{(\Delta+i\omega)^2+ \omega_0^2}\right),
\]
and therefore $\text{Im}(D_{\frac K2 \widehat g_{\Delta, \omega_0}} (\omega))=0$ iff $\omega=0$ or $\omega= \pm \sqrt{\omega_0^2-\Delta^2}$. The condition \eqref{PENROSE} then reads $K<K_{\Delta,\omega_0}$ where the threshold $K_{\Delta,\omega_0}$ is given by
\begin{equation}
\label{eq:KcDelta}
K_{\Delta, \omega_0}\, =\, \begin{cases}
\frac{2 (\Delta^2+ \omega_0^2)}\Delta = \frac 2{\pi g_{\Delta, \omega_0}(0)} & \text{ if } \omega_0\leq \Delta \, , \\
4 \Delta &\text{ if } \omega_0 > \Delta\,.
\end{cases}
\end{equation}
In addition, stationary or periodic solutions with $|R|\neq 0$ exist for $K>K_{\Delta,\omega_0}$ \cite{MBSOSA09}.

Finally, notice that we do not know whether conditions \eqref{STABCRIT} and \eqref{PENROSE} are equivalent in all cases. However, their equivalence is not limited to symmetric examples above and one can show that it holds for every density $\alpha g_\Delta(\cdot+\omega_0)+(1-\alpha)g_\Delta(\cdot-\omega_0)$ where $\alpha\in [0,1]$ is arbitrary.

\section{Dynamics of the (rescaled) order parameter}
This section aims at establishing a Volterra integral equation for a rescaled order parameter and to use this equation in order to obtain an estimate for the quantity $\sup\limits_{t\in [0,T]}(1+t)^n|R(t)|$, independently of $T>0$. The first part of the procedure is similar to the one in \cite{SMM92}. The approach here follows even more closely the Vlasov equation analysis as it also employs the Galilean transformation $(t,\theta,\omega)\mapsto (t,\theta+t\omega,\omega)$ prior to the Fourier transform. The second part reproduces in details the methods presented in \cite{FR14,V10}.

\subsection{Volterra equation for the (rescaled) order parameter}
Instead of the decomposition \eqref{DECOMPER}, it turns out more convenient to separate the perturbation norms from the functions themselves (no matter which norm is involved) \cite{SM91}, namely we write 
\[
\rho(t,\theta,\omega)=\frac{1}{2\pi}+\epsilon r(t,\theta,\omega),\ \forall t\in\R^+.
\]
where $\epsilon>0$ (and the norm of $r$ ought to be prescribed). In other words, we have $\text{r}=\epsilon r$ and $\text{R}=\epsilon R$ where the rescaled order parameter $R$ is defined by 
\[
R(t)= \int_{\T^1\times\R} r(t,\theta,\omega)g(\omega)e^{-i\theta}d\theta d\omega, \forall t\in\R^+.
\] 
Plugging the ansatz $\rho=\frac{1}{2\pi}+\epsilon r$ into equation \eqref{BASICDYNAM} and applying the mentioned Galilean transformation, the following PDE results for the quantity $p(t)=T^tr(t)\cdot g$ explicitly given by
\[
p(t,\theta,\omega)=r(t,\theta+t\omega,\omega)g(\omega),\ \forall (\theta,\omega)\in \T^1\times\R,\ t\in\R^+,
\]
\begin{equation}
\partial_tp(\theta,\omega)+\epsilon\partial_\theta p(\theta,\omega)W(\theta+t\omega,p)+\left(\frac{g(\omega)}{2\pi}+\epsilon p(\theta,\omega)\right)\partial_\theta W(\theta+t\omega,p)=0,\ \forall (\theta,\omega)\in \T^1\times\R,\ t>0,
\label{GALILEANVAR}
\end{equation}
where
\[
W(\theta,p)=K\int_{\T^1\times\R}\sin (\vartheta+t\omega-\theta)p(\vartheta,\omega)d\vartheta d\omega.
\]

Now, given any $L^1$-function $u$ of the cylinder, let $\widehat{u}_k(\tau)$ denote its Fourier transform defined by
\[
\widehat{u}_k(\tau)=\int_{\T^1\times\R}u(\theta,\omega)e^{-i(k\theta+\tau\omega)}d\theta d\omega,\ \forall (k,\tau)\in \Z\times\R.
\]
(We obviously have $\widehat{g}=\frac1{2\pi}\widehat{(1_{\T^1}\cdot g)}_0$ when the product $1_{\T^1}\cdot g$ is regarded as a function defined on the cylinder.)

The solution $p(t)$ of equation \eqref{GALILEANVAR} must be absolutely integrable over the cylinder for every $t\in\R^+$. Observing that $R(t)=\widehat{p}_1(t,t)$, from \eqref{GALILEANVAR}, one can derive the following infinite system of coupled ODEs for the quantities $\{\widehat{p}_k(t,\tau)\}$ (where we include the explicit dependence on time for clarity)
\begin{equation}
\partial_t \widehat{p}_k(t,\tau)+\frac{kK}{2}\left(\overline{R(t)}\left(\widehat{g}(\tau+t)\delta_{k,-1}+\epsilon\widehat{p}_{k+1}(t,\tau+t)\right)-R(t)\left(\widehat{g}(\tau-t)\delta_{k,1}+\epsilon\widehat{p}_{k-1}(t,\tau-t)\right)\right)=0,
\label{FOURIERVAR}
\end{equation}
for all $(k,\tau)\in\Z\times\R$ and $t>0$, where we have used the Kronecker symbol. (NB: For $k=0$, we always have $\widehat{p}_0(t,\tau)=0$, as a consequence of the constraint $\int_{\T^1}r(t,\theta,\omega) d\theta =0$ for all $\omega\in\R$.)
Letting $k=1$, integrating in time and letting $\tau=t$, we finally obtain the desired equation for the rescaled order parameter
\begin{equation}
R(t)-\frac{K}{2}\int_0^t\widehat{g}(t-s)R(s)ds=F(t),\forall t\in\R^+,\ \text{where}\ F(t)=\widehat{p}_1(0,t)-\frac{\epsilon K}{2}\int_0^t\widehat{p}_2(s,t+s)\overline{R(s)}ds.
\label{FULLEQ}
\end{equation}
Regarding the term $F$ as an autonomous input signal, equation \eqref{FULLEQ} appears to be a Volterra equation of the second kind. As such, it is well-known to have a unique solution for every density $g\in L^1(\R)$ \cite{V59}. More importantly for our purpose, tricky arguments based on Complex Analysis, and inspired by \cite{FR14,V10}, show that under suitable stability conditions (as listed in Theorem \ref{MAINRES}), the solution polynomial decay is controlled by the input. As exposed in the next section, this control holds in a broader context than the strict analysis of the Kuramoto model.

\subsection{Polynomial decay of solutions of Volterra equations. Application to equation \eqref{FULLEQ}}
Recall that $\Pi^-$ denotes the lower half plane of complex number $z$ with $\text{Im}(z)\leq 0$. For any $G\in L^1(\R^+)$, the function $D_G$ defined by \eqref{DEFDG} can be extended to a well-defined function on $\Pi^-$ (which we denote by the same symbol).
\begin{Pro}
Let $n\in\N$ and let $G\in L^1(\R^+)\cap L^\infty(\R^+)$ be a complex valued function that satisfies 
\[
\int_{\R^+}t^4|G(t)|^2dt<+\infty,\quad \int_{\R^+}t^n|G(t)|dt<+\infty\quad\text{and}\quad D_G|_{\Pi^-}\neq 0.
\]
There exists $C_{n,G}\in\R^+$ such that for every complex valued function $F$ on $\R^+$, the solution of the Volterra equation
\begin{equation}
R(t)=F(t)+\int_0^tG(t-s)R(s)ds,\ \forall t\in\R^+,
\label{VOLTERRA}
\end{equation}
satisfies the following inequality, for any $T > 0$ :
\begin{equation}
\sup_{t\in [0,T]}(1+t)^n|R(t)|\leq C_{n,G}\sup_{t\in [0,T]}(1+t)^n|F(t)|,\ \forall T>0.
\label{FINITIM}
\end{equation}
\label{GENCONTROL}
\end{Pro}
The proof is given in section \ref{S-GENCONTROL} below. 
The following comments prepare its application to equation \eqref{FULLEQ} of the Kuramoto model. 

For $n\geq 4$, the condition $\int_{\R^+}t^4|G(t)|^2dt<+\infty$ is redundant as it can be deduced from the assumptions $G\in L^1(\R^+)\cap L^\infty(\R^+)$ and $\int_{\R^+}t^n|G(t)|dt<+\infty$.

Moreover, we have $G=\frac{K}2\widehat{g}$ in equation \eqref{FULLEQ}; hence it suffices to impose $\widehat{g}\in L^1(\R^+)$, $\int_{\R^+}t^4|\widehat{g}(t)|^2dt<+\infty$, $\int_{\R^+}t^n|\widehat{g}(t)|dt<+\infty$ and condition \eqref{STABCRIT}, in order to apply Proposition \ref{GENCONTROL} there. Accordingly, the following statement immediately results, which we only state for $n\geq 4$ anticipating the condition that will result from the bootstrap argument below in the proof of Theorem \ref{MAINRES}. In addition, anticipating also the proof of Proposition \ref{LARGEPERT}, we explicitly express the dependence on $K$ in the constant $C_K$ of the statement (which readily follows from the expression of $C_{n,G}$ in the end of the proof of Proposition \ref{GENCONTROL}).
\begin{Cor}
Assume that $g$ is of class $C^n(\R)$ for some $n\geq 4$ and satisfies 
\[
\widehat{g}\in L^1(\R^+)\quad\text{and}\quad\int_{\R^+}t^n|\widehat{g}(t)|dt<+\infty.
\]
For every $K\geq 0$ such that the condition \eqref{STABCRIT} holds,  
there exists $C_K>0$ such that for every input signal $F$, the solution of equation \eqref{FULLEQ} satisfies the following property
\[ 
\sup_{t\in [0,T]}(1+t)^n|R(t)|\leq C_K\sup_{t\in [0,T]}(1+t)^n|F(t)|,\ \forall T>0.
\]
The constant $C_K$ is a polynomial of order $n+1$ in $K$ with $g$-dependent coefficients.
\label{FIRSTINEQ}
\end{Cor}

\subsection{Proof of Proposition \ref{GENCONTROL}}\label{S-GENCONTROL}
Following \cite{FR14}, the proof of Proposition \ref{GENCONTROL} proceeds by iterations on the integer $n$ and thus starts with $n=0$. This first step itself separates into two parts; first, we show (in the next statement below)  that the solution $R$  of the Volterra equation is square integrable if the source term $F$ is. Then, we  use this square integrability to establish the estimate of Proposition \ref{GENCONTROL}, by application of the Fourier transform.
\subsubsection{Proof of the estimate for $n=0$}
\begin{Pro}
Let $G\in L^1(\R^+)\cap L^\infty(\R^+)$ be such that
\[
\int_{\R^+}t^4|G(t)|^2dt<+\infty\quad\text{and}\quad D_G|_{\Pi^-}\neq 0.
\]
If the input signal $F\in L^2(\R^+) \cap L^\infty(\R^+)$, then the solution $R$ of equation \eqref{VOLTERRA} also belongs to $L^2(\R^+) \cap L^\infty(\R^+)$.
\label{PRELIMBOUNDS}
\end{Pro}

\noindent
{\sl Proof}. We shall actually prove a slightly refined result, namely:

\noindent
(i) {\em  There exists $C_G$ such that for any $F \in L^2(\R^+)$, we have $R \in L^2(\R^+)$ with estimate  
\[
\|R\|_{L^2(\R^+)}\leq C_G \|F\|_{L^2(\R^+)}.
\]}

\noindent
(ii) {\em There exists $C'_G > 0$ such that for any $F \in L^2(\R^+) \cap L^\infty(\R^+)$, we have $R \in L^2(\R^+) \cap L^\infty(\R^+)$ with estimate
\[
\|R\|_{L^\infty(\R^+)}\leq \|F\|_{L^\infty(\R^+)} + C'_G \|F\|_{L^2(\R^+)}.
\]}
Most of the proof consists in showing statement (i). Indeed, once this is proved, statement (ii) will follow immediately with $C'_G=C_G\|G\|_{L^2(\R^+)}$, after applying the 
Cauchy-Schwarz inequality to the integral term in the right hand side of equation \eqref{VOLTERRA}. 

For the proof of (i), we need to introduce the trivial extension $\widetilde{u}$ to $\R$ of a function on $\R^+$, defined as follows
\[
\widetilde{u}|_{\R^+}=u\quad\text{and}\quad \widetilde{u}_{\R^-_\ast}=0.
\]
Now observe that for every solution $R\in L^2(\R^+)$ of equation \eqref{VOLTERRA}, $\widetilde{R}$ satisfies
\begin{equation}
\widetilde{R}=\widetilde{F}+\widetilde{G}\ast\widetilde{R},
\label{FOURVOLT}
\end{equation}
over $\R$,  and its Fourier transform can be easily solved by the expression $\frac{\widehat{\widetilde F}}{1-\widehat{\widetilde G}}$, provided that $(1-\widehat{\widetilde G})|_{\R}\neq 0$, {\sl viz.}\ $D_G|_{\R}\neq 0$ since $1-\widehat{\widetilde G}=D_G$. However, the condition $G\in L^1(\R^+)$ implies that $\widehat{\widetilde G}$ is continuous on $\R$ and $\widehat{\widetilde G}(\pm\infty)=0$. Hence, the condition $D_G|_{\R}\neq 0$ implies the existence of $\delta>0$ such that $D_G|_{\R}\geq \delta$ and the function $\frac{\widehat{\widetilde F}}{D_G}$ is in $L^2(\R)$ for every $F\in L^2(\R^+)$.

\noindent
By the Plancherel Theorem, the function $R_0\in L^2(\R)$ defined via the inverse Fourier transform ${\cal F}^{-1}$ as follows
\[
R_0={\cal F}^{-1}\left(\frac{\widehat{\widetilde F}}{D_G}\right),
\]
satisfies $\|R_0\|_{L^2(\R)}\leq \frac{\sqrt{2\pi}}{\delta}\|F\|_{L^2(\R^+)}$. We are going to prove that $\|R_0\|_{L^2(\R^-)}=0$. Uniqueness of the solutions of equation \eqref{VOLTERRA} will then imply $\|R\|_{L^2(\R^+)}=\|R_0\|_{L^2(\R)}$ and the inequality in statement (i) will hold with $C_G=\frac{\sqrt{2\pi}}{\delta}$. 

To that goal, we shall need the following considerations. First, the easy part of the Paley-Wiener Theorem implies that the function $\widehat{\widetilde F}$ extends to an analytic function on the lower half plane 
\[
\Pi^-_\ast=\{x-iy\ :\ x\in\R, y\in\R^+_*\},
\]
namely,
\begin{equation} \label{paleywiener}
\widehat{\widetilde F}(x-iy)=\int_{\R^+}F(t)e^{-yt}e^{-ixt}dt = \widehat{( \widetilde{F}e^{-y \cdot})}(x).
\end{equation}
 Second, we have the following statement. 
\begin{Lem}
For every $F\in L^2(\R^+)$, we have 
\[
\lim_{y\to +\infty}\sup_{x\in\R}|\widehat{\widetilde F}(x-iy)|=0\quad\text{and}\quad
\lim_{x\to \pm\infty}\sup_{y>\epsilon}|\widehat{\widetilde F}(x-iy)|=0,\ \forall \epsilon>0.
\]
Moreover, if we also have $F\in L^1(\R^+)$, then the former property also holds for $\epsilon=0$.
\label{ASYMPBEHA}
\end{Lem}
{\sl Proof of the Lemma.} The first limit easily follows by applying the Cauchy-Schwarz inequality to the integral in the right hand side of \eqref{paleywiener}. 
For the second limit, it suffices to apply the uniform Riemann-Lebesgue Lemma, see {\sl e.g.}\ \cite{P02}, to the family $\{F_y\}_{y\in [\epsilon,+\infty]}$ of functions defined by 
\[
F_y(t)=\left\{\begin{array}{cl}
F(t)e^{-yt}&\text{if}\ y\in [\epsilon,+\infty)\\
0&\text{if}\ y=+\infty
\end{array}\right.\ \forall t\in\R^+.
\]
Each function $F_y\in L^1(\R^+)$ (this can be seen by applying Cauchy-Schwarz and using $F\in L^2(\R^+)$). All we have to show is that the family is compact in $L^1(\R)$. Since $[\epsilon,+\infty]$ is compact after addition of the point $+\infty$, it suffices to show that $y\mapsto F_y$ is continuous in $L^1(\R^+)$. This is granted by the following expression
\[
\int_{\R^+}|F_y(t)-F_{y'}(t)|dt\leq \|F\|_{L^2(\R^+)}\sqrt{\int_{\R^+}(e^{-yt}-e^{-y't})^2dt}=\|F\|_{L^2(\R^+)}\sqrt{\frac1{2y}+\frac1{2y'}-\frac2{y+y'}}.
\]
\hfill $\Box$

Since $G\in L^2(\R^+)$, Lemma \ref{ASYMPBEHA} together with the condition $D_G|_{\Pi^-}\neq 0$ implies the existence of $\delta>0$ such that $D_G|_{\Pi^-_\ast}\geq \delta$. 
Invoking again the Paley-Wiener Theorem, the function $\frac{\widehat{\widetilde F}}{D_G}$ is analytic on the lower half plane $\Pi^-_\ast$. Consider the function $R_\epsilon\in L^2(\R^-)$ defined for $\epsilon>0$ by 
\[
R_\epsilon(t)=\frac1{2\pi}\int_{\R}\frac{\widehat{\widetilde F}(x-i\epsilon)}{D_G(x-i\epsilon)} e^{ixt}dx,\ \text{for a.e.}\  t\in\R^-.
\]
\begin{Lem} 
{\rm (a)} $\lim_{\epsilon\searrow 0}\|R_\epsilon-R_0\|_{L^2(\R^-)}=0$ and 
{\rm (b)}
$\|R_\epsilon\|_{L^2(\R^-)}=0$ for every $\epsilon>0$. 
\label{TECHNIK}
\end{Lem}
Lemma~\ref{TECHNIK} directly yields $\|R_0\|_{L^2(\R^-)}=0$.  This completes the proof of statement (i) of Proposition \ref{PRELIMBOUNDS} and hence the proof of the proposition is complete. \hfill $\Box$

\noindent
{\sl Proof of  Lemma~\ref{TECHNIK}.} (a) One easily obtains the following inequality
\begin{equation*}
\left\|\frac{\widehat{\widetilde F}(\cdot-i\epsilon)}{D_G(\cdot-i\epsilon)}-\frac{\widehat{\widetilde F}(\cdot)}{D_G(\cdot)}\right\|_{L^2(\R)} 
\leq \frac{1}{\delta} \left\|  \widehat{\widetilde F}(\cdot-i\epsilon) - \widehat{\widetilde F} (\cdot) \right\|_{L^2(\R)} + \frac{1}{\delta^2} \left\|  \widehat{\widetilde F} (\cdot)\, \left( \widehat{\widetilde G}(\cdot-i\epsilon) - \widehat{\widetilde G}(\cdot) \right)\right\|_{L^2(\R)} ,
\end{equation*}
The first term vanishes when $\epsilon\to 0$. Indeed, by the Plancherel formula, we have
\[
\left\|  \widehat{\widetilde F}(\cdot-i\epsilon) - \widehat{\widetilde F} (\cdot) \right\|_{L^2(\R)}
= \sqrt{2\pi} \, \| \widetilde{F} (\cdot) e^{-\eps \cdot} - \widetilde F (\cdot) \|_{L^2(\R)} ,
\]
and the right-hand side vanishes by the dominated convergence theorem. The second term also vanishes by combining this same theorem with the continuity and boundedness of $\widehat{\widetilde G}$. It results that
\[
\lim_{\epsilon\to 0}\left\|\frac{\widehat{\widetilde F}(\cdot-i\epsilon)}{D_G(\cdot-i\epsilon)}-
\frac{\widehat{\widetilde F} (\cdot)}{D_G (\cdot)}\right\|_{L^2(\R)}=0.
\]
The result then immediately follows from the continuity of the Fourier transform in $L^2$ topology.

\noindent
(b) It is equivalent to prove that
\[
R(t) e^{\eps t} = \frac1{2\pi}\int_{\R}\frac{\widehat{\widetilde F}(x-i\epsilon)}{D_G(x-i\epsilon)}e^{ixt+\eps t}dx,
\]
vanishes for almost every $t$ in $\R^-$.

The function $z \mapsto \frac{\widehat{\widetilde F}(z)}{D_G(z)}e^{izt}$ is analytic. Hence, given $s>0$ and $\mu>\epsilon$, the Cauchy theorem implies that the integral $\int_{-s}^s \frac{\widehat{\widetilde F}(x-i\epsilon)}{D_G(x-i\epsilon)}e^{ixt+\epsilon t}dx$ is equal to the following sum
\[
-i e^{-ist}\int_\epsilon^\mu\frac{\widehat{\widetilde F}(-s-iy)}{D_G(-s-iy)}e^{y t}dy 
+i e^{ist}\int_\epsilon^\mu\frac{\widehat{\widetilde F}(s-iy)}{D_G(s-iy)}e^{y t}dy+
e^{\mu t}\int_{-s}^s\frac{\widehat{\widetilde F}(x-i\mu)}{D_G(x-i\mu)}e^{ixt}dx.
\]
All we have to do next is to show that each term in this decomposition vanishes when taking the limit $\mu\to +\infty$ and then $s\to +\infty$. Since we have $D_G|_{\Pi^-}\geq \delta$, we can forget about denominators and focus on the behavior of the numerators in this process.

Now, we have
\[
\lim_{\mu\to+\infty} \widehat{\widetilde F}(x-i\mu)=0,\ \forall x\in \R,
\]
hence the third term vanishes in the limit. Moreover, one can replace $\mu$ by $+\infty$ in the first two terms, using $t<0$ and the fact that the modulus $|\widehat{\widetilde F}(\pm s-iy)|$ remains bounded. Finally, that the first two terms vanish in the limit $s\to+\infty$ is a consequence of the second limit behavior in Lemma \ref{ASYMPBEHA} (together with the dominated convergence theorem). The proof of Lemma~\ref{TECHNIK} is complete. \hfill $\Box$

\noindent
{\sl Proof of Proposition \ref{GENCONTROL}}. We start with the case $n=0$. Let $\chi:\R\to [0,1]$ be a smooth function such that 
\[
\chi(x)=\left\{\begin{array}{cl}
1&\text{if}\ |x|\leq \tfrac12\\
0&\text{if}\ |x|\geq 1
\end{array}\right.
\]
and given $\eta>0$, let $\chi_\eta$ be the function defined by 
\[
\chi_\eta(x)=\chi\left(\tfrac{x}{\eta}\right),\ \forall x\in\R.
\]
Let $T_\eta$ be the operator defined in $L^\infty(\R)$ by the convolution with kernel given by the inverse Fourier transform ${\cal F}^{-1}\chi_\eta$.

We are going to prove the existence of $C_{0,G}>0$ such that for every $F\in L^2(\R)\cap L^\infty(\R)$, the solution of equation \eqref{VOLTERRA} satisfies the following inequality 
\begin{equation}
\|R\|_{L^\infty(\R^+)}\leq C_{0,G}\|F\|_{L^\infty(\R^+)}.
\label{ANOTINEQ}
\end{equation}
Proposition \ref{GENCONTROL} for $n=0$ will then follow by applying this inequality to the input signal $F_T = F 1_{[0,T]}$ (here, $T>0$ is arbitrary and $1_{[0,T]}$ denotes the characteristic function of $[0,T]$) and associated solution $R_T$, thanks to the fact that $F_T\in L^2(\R)\cap L^\infty(\R)$ {\sl for every} input signal $F$. Indeed, using both that $F_T=F$ on $[0,T]$ and uniqueness of solutions of Volterra equation, we obtain $R_T=R$ on $[0,T]$ and then
\[
\|R\|_{L^\infty([0,T])}= \|R_T\|_{L^\infty([0,T])}\leq \|R_T\|_{L^\infty(\R^+)}\leq C_{0,G}\|F_T\|_{L^\infty(\R^+)}=C_{0,G}\|F\|_{L^\infty([0,T])},
\]
as desired.

From now on, consider an input signal $F\in L^2(\R)\cap L^\infty(\R)$ and recall that the trivial extension $\widetilde R$ of the solution of the Volterra equation \eqref{VOLTERRA}. By Proposition \ref{PRELIMBOUNDS}, we have $\widetilde R\in L^2(\R)\cap L^\infty(\R)$; hence the following functions are well-defined
\[
\check{R}_H=(\text{Id}-T_\eta)\widetilde R\quad\text{and}\quad \check{R}_L=T_\eta\widetilde R,
\]
and we seek for bounds on their uniform norm that are independent of $\|F\|_{L^2(\R^+)}$. 

Applying the operator $(\text{Id}-T_\eta)$ to equation \eqref{FOURVOLT}, and then Fourier transform, we obtain using also the relation $1-\chi_\eta=(1-\chi_\eta)(1-\chi_{\frac{\eta}2})$
\[
\widehat{\check{R}}_H=(1-\chi_\eta)\widehat{\widetilde F}+(1-\chi_{\frac{\eta}2})\widehat{\widetilde G}\widehat{\check{R}}_H,
\]
from where the following relation immediately follows
\[
\check{R}_H=(\text{Id}-T_\eta)\widetilde F+(\text{Id}-T_{\frac{\eta}2})\widetilde G\ast \check{R}_H.
\]
The assumption on $\chi$ implies $\|{\cal F}^{-1}\chi_\eta\|_{L^1(\R)}<+\infty$. Using $\|T_\eta\widetilde F\|_{L^\infty(\R)}\leq \|{\cal F}^{-1}\chi_\eta\|_{L^1(\R)}\|\widetilde F\|_{L^\infty(\R)}$, the following inequality results
\[
(1-\|(\text{Id}-T_{\frac{\eta}2})\widetilde G\|_{L^1(\R)})\|\check{R}_H\|_{L^\infty(\R)}\leq (1+\|{\cal F}^{-1}\chi_\eta\|_{L^1(\R)})\|\widetilde F\|_{L^\infty(\R)}.
\]
We are going to prove that $\lim\limits_{\eta\to +\infty}\|(\text{Id}-T_\eta)\widetilde G\|_{L^1(\R)}=0$. This implies 
\begin{equation}
\|\check{R}_H\|_{L^\infty(\R)}\leq 2(1+\|{\cal F}^{-1}\chi_\eta\|_{L^1(\R)})\|\widetilde F\|_{L^\infty(\R)},
\label{ESTIMHIGH}
\end{equation}
provided that $\eta$ is sufficiently large. 

Recall the notation $\langle t\rangle =\sqrt{1+t^2}$ and let $Du$ denotes the derivative of the function $u$. By the Plancherel Theorem, we have ($\|\langle t\rangle^2 (\text{Id}-T_\eta)\widetilde G\|_{L^2(\R)}$ is finite thanks to the assumptions on $G$)
\begin{multline*}
\frac{\|\langle t\rangle^2 (\text{Id}-T_\eta)\widetilde G\|_{L^2(\R)}}{\sqrt{2\pi}}=\\
\left\|(1-D^2)\left((1-\chi_\eta)\widehat{\widetilde G}\right)\right\|_{L^2(\R)}\leq \left\|(1-\chi_\eta)\widehat{\widetilde G}\right\|_{L^2(\R)}+\left\|D^2\left((1-\chi_\eta)\widehat{\widetilde G}\right)\right\|_{L^2(\R)}.
\end{multline*}
The first term in the sum is bounded above by 
\[
\int_{\R\setminus [-\frac{\eta}2,\frac{\eta}2]}\Big|\widehat{\widetilde G}(t)\Big|^2dt,
\] 
which vanishes as $\eta\to +\infty$. The second term has to be finite thanks to the assumption on $G$ and the fact that $\|D^j\chi_\eta\|_{L^\infty(\R)}=O(\eta^{-j})$ for $j\in \{1,2\}$. This asymptotic behavior implies that the contributions involving the derivatives $D^j\chi_\eta$ must vanish as $\eta\to +\infty$. For the remaining contribution $(1-\chi_\eta)D^2\widehat{\widetilde G}$, the assumption $\int_{\R^+}t^4|G(t)|^2dt<+\infty$ implies
\[
\lim_{\eta\to +\infty}\left\Vert(1-\chi_\eta)D^2\widehat{\widetilde G}\right\Vert_{L^2(\R)}\leq \lim_{\eta\to +\infty}\int_{\R\setminus [-\frac{\eta}2,\frac{\eta}2]}\left|D^2\widehat{\widetilde G}(t)\right|^2dt=0,
\]
from where it follows that the second term in the sum above also vanishes in the limit $\eta\to +\infty$. The desired behavior then follows from the inequality
\[
\|(\text{Id}-T_\eta)\widetilde G\|_{L^1(\R)}\leq \left(\int_{\R}\frac{dt}{\langle t\rangle^4}\right)^{\frac12}\|\langle t\rangle^2(\text{Id}-T_\eta)\widetilde G\|_{L^2(\R)}.
\]
In order to estimate $\|\check{R}_L\|_{L^\infty(\R)}$, we observe that the expression $\widehat{\widetilde R}=\frac{\widehat{\widetilde F}}{D_G}$ yields
\[
\widehat{\check R}_L=\frac{\chi_\eta}{D_G}\widehat{\widetilde F},
\]
and then
\[
{\check R}_L={\cal F}^{-1}\left(\frac{\chi_\eta}{D_G}\right)\ast \widetilde F,
\]
which implies the existence of $C_{\eta,G}>0$ such that the following inequality holds
\[
\|\check{R}_L\|_{L^\infty(\R)}\leq C_{\eta,G} \|\widetilde F\|_{L^\infty(\R)}.
\]
Adding this estimate with the one in equation \eqref{ESTIMHIGH}, the inequality \eqref{ANOTINEQ} follows with 
\[
C_{0,G}=2(1+\|{\cal F}^{-1}\chi_\eta\|_{L^1(\R)})+C_{\eta,G},
\]
and Proposition \ref{GENCONTROL} is proved for $n=0$.

\subsubsection{Proof of the estimate for $n=1$}
In order to prove Proposition \ref{GENCONTROL} for $n=1$, we observe that the function $R_1$ defined by $R_1(t)=(1+t)R(t)$ satisfies the equation
\[
R_1(t)=F_1(t)+\int_0^tG(t-s)R_1(s)ds,\ \forall t\in\R^+,
\]
where 
\[
F_1(t)=(1+t)F(t)+(\widetilde{G}_1\ast\widetilde R)(t)\quad\text{and}\quad G_1(t)=t G(t).
\]
Therefore, it suffices to check that $\widetilde{G}_1\ast\widetilde R\in L^2(\R^+)\cap L^\infty(\R^+)$ in order to apply the result for $n=0$ with input term $F_1$. (Recall from the proof for $n=0$ that we may assume without loss of generality that $(1+t)F\in L^2(\R^+)\cap L^\infty(\R^+)$.)
The additional assumption in the statement of the proposition when passing from $n=0$ to $n=1$ is actually $G_1\in L^1(\R^+)$. Applying adequately the Young inequality, we obtain
\[
\|\widetilde{G}_1\ast\widetilde R\|_{L^2(\R^+)}\leq \|G_1\|_{L^1(\R^+)}\|R\|_{L^2(\R^+)}\quad\text{and}\quad \|\widetilde{G}_1\ast\widetilde R\|_{L^\infty(\R^+)}\leq \|G_1\|_{L^1(\R^+)}\|R\|_{L^\infty(\R^+)},
\]
from where the desired conclusion follows from $R\in L^2(\R^+)\cap L^\infty(\R^+)$ (consequence of Proposition \ref{PRELIMBOUNDS}). Therefore Proposition \ref{GENCONTROL} holds for $n=1$ with 
\[
C_{1,G}=C_{0,G}(1+C_{0,G}\|G_1\|_{L^1(\R^+)}).
\]

\subsubsection{Proof of the estimate for $n>1$}
For $n>1$, the argument is similar. The function $R_n(t)=(1+t)^nR(t)$ satisfies the equation
\[
R_n(t)=F_n(t)+\int_0^tG(t-s)R_n(s)ds,\ \forall t\in\R^+,
\]
where 
\[
F_n(t)=(1+t)^nF(t)+\int_0^t\left((1+t)^n-(1+s)^n\right)G(t-s)R(s)ds.
\]
Using the inequality $(1+t)^n-(1+s)^n\leq (2^n-1)\left((t-s)(1+s)^{n-1}+(t-s)^n\right)$, we obtain
\begin{multline*}
\left|\int_0^t\left((1+t)^n-(1+s)^n\right)G(t-s)R(s)ds\right|\leq 
\\
(2^n-1)\left(\int_0^t|G_1(t-s)||R_{n-1}(s)|ds+\int_0^t|G_n(t-s)||R(s)|ds\right),
\end{multline*}
where $G_n(t)=t^nG(t)$. Similarly to as for the proof in the case $n=1$, it therefore suffices to ensure
\[
\widetilde{G}_1\ast \widetilde{R}_{n-1} \in L^2(\R^+)\cap L^\infty(\R^+)\quad\text{and}\quad
\widetilde{G}_n\ast \widetilde{R} \in L^2(\R^+)\cap L^\infty(\R^+).
\]
The first property follows from an induction argument; hence the (additional) conditions $G_k\in L^1(\R^+)$ for $k\in \{2,\cdots ,n-1\}$. The second property is a consequence of the assumption $G_n\in L^1(\R^+)$. We conclude that the statement holds with 
\[
C_{n,G}=C_{0,G}\left(1+(2^n-1)(C_{n-1,G}\|G_1\|_{L^1(\R^+)}+C_{0,G}\|G_n\|_{L^1(\R^+)})\right),
\]
and the proof of Proposition \ref{GENCONTROL} is complete.
\hfill $\Box$

\subsection{Additional results}
\subsubsection{Simple proof of Proposition \ref{GENCONTROL} under stronger constraint}
By applying adequately the Young inequality to equation \eqref{FOURVOLT}, one gets
\[
\|R\|_{L^\infty(\R^+)}\leq \|F\|_{L^\infty(\R^+)}+\|G\|_{L^1(\R^+)}\|R\|_{L^\infty(\R^+)}.
\]
Accordingly, the conclusion of Proposition \ref{GENCONTROL} for $n=0$  follows by assuming $\|G\|_{L^1(\R^+)}<1$ instead of 
\[
\int_{\R^+}t^4|G(t)|^2dt<+\infty,\quad \text{and}\quad D_G|_{\Pi^-}\neq 0.
\]
(The Proposition itself then follows when assuming the other conditions $G\in L^1(\R^+)\cap L^\infty(\R^+)$ and $\int_{\R^+}t^n|G(t)|dt<+\infty$.) The condition $\|G\|_{L^1(\R^+)}<1$ is in general stronger than $D_G|_{\Pi^-}\neq 0$. However, notice that when $g$ is unimodal and such that $\widehat g \ge 0$, both conditions $\|\hat{g}\|_{L^1(\R^+)}<1$ and $D_{\hat{g}}|_{\Pi^-}\neq 0$ are equivalent. This applies in particular to Gaussian and Cauchy distribution densities. 

\subsubsection{Optimality of the stability criterion}\label{S-OPTIMAL}
As announced in section \ref{S-ANALSTABCOND}, one can show exponential growth of solutions of the Volterra equation when the function $D_G$ associated with the kernel $G$, vanishes at some point in the lower half plane $\Pi^-_\ast$.
\begin{Lem}
Let $G\in L^2(\R^+)$ be such that $D_G(\omega_0)=0$ for some $\omega_0\in \Pi^-_\ast$. Then, there exists $F\in L^2(\R^+)\cap L^\infty(\R^+)$ such that the solution of the Volterra equation \eqref{VOLTERRA} writes
\[
R(t)=Ae^{i\omega_0t},\ \forall t_0\in\R^+,
\]
(notice that $\text{Re}(i\omega_0)>0$) where $A\in\R$, $A\neq 0$. 
\label{OPTIMAL}
\end{Lem}
{\sl Proof.} Given $A\in\R$ and $\omega\in\C$, the function $t\mapsto Ae^{i\omega t}$ solves the Volterra equation iff we have 
\[
F(t)=Ae^{i\omega t}\left(1-\int_0^tG(s)e^{-i\omega s}ds\right),\ \forall t\in\R^+.
\]
For $\omega=\omega_0$, by using $D_G(\omega_0)=0$, this expression simplifies to the following one
\[
F(t)=A\int_{\R^+}G(t+s)e^{-i\omega_0 s}ds,\ \forall t\in\R^+.
\]
By the Cauchy-Schwarz inequality
\[
|F(t)|\leq A\|G\|_{L^2(\R^+)}\|e^{\text{Im}(\omega_0\cdot)}\|_{L^2(\R^+)},\ \forall t\in\R^+,
\]
which implies that $F\in L^\infty(\R^+)$. Moreover, one obtains by  Fubini's Theorem 
\begin{align*}
\|F\|_{L^2(\R^+)}^2&\leq A^2\int_{\R^+\times\R^+\times\R^+}|G(t+s_1)||G(t+s_2)|e^{\text{Im}(\omega_0) s_1}e^{\text{Im}(\omega_0) s_2}ds_1ds_2dt
\\
&\leq A^2\int_{\R^+\times\R^+}\left(\int_{\R^+}|G(t+s_1)||G(t+s_2)|dt\right)e^{\text{Im}(\omega_0) s_1}e^{\text{Im}(\omega_0) s_2}ds_1ds_2\\
&\leq A^2 \|G\|_{L^2(\R^+)}^2\|e^{\text{Im}(\omega_0\cdot)}\|_{L^1(\R^+)}^2.
\end{align*}
\hfill $\Box$

\subsubsection{Proof of Lemma \ref{EQUISTAB}}\label{S-EQUISTAB}
We have $D_G|_{\R}\neq 0$, hence by applying Lemma \ref{ASYMPBEHA} to $G$, we conclude that every zero of $D_G$ in $\Pi^-$ must lie in the interior of a bounded rectangle 
\[
\left\{x-iy:\, |x|\leq s_G, 0\leq y\leq \mu_G\right\},
\]
for some suitable $s_G,\mu_G>0$.  Moreover, the function $D_G$ is analytic in $\Pi^-_\ast$. By the argument principle, given $\epsilon>0$, the number $N_\epsilon$ of zeros of $D_G$ in the  half plane
$
\{x-iy:\, x\in\R,\,  y>\epsilon\}$
is given by
\[
\frac1{2\pi i}\left(\int_{-s}^s\frac{-\widehat{\widetilde G}'(x-i\epsilon)}{D_G(x-i\epsilon)}dx+i\int_{\epsilon}^\mu\frac{-\widehat{\widetilde G}'(s-iy)}{D_G(s-iy)}dy-\int_{-s}^s\frac{-\widehat{\widetilde G}'(x-i\mu)}{D_G(x-i\mu)}dx\right.
\left. -i\int_{\epsilon}^\mu\frac{-\widehat{\widetilde G}'(-s-iy)}{D_G(-s-iy)}dy\right),
\]
for every $s\geq s_G$ and $\mu\geq \mu_G$. As in the proof of Lemma \ref{TECHNIK}, the third term vanishes when $\mu\to +\infty$. For the second term, we remark that  $i\widehat{\widetilde G}'$ is the Fourier transform of $t\mapsto tG(t)$. It follows  easily that 
\[
 \left| \frac{-\widehat{\widetilde G}'(s-iy)}{D_G(s-iy)} \right| \: \leq \: \frac 1 c \int_{\R^+} | t G(t)| e^{-yt} dt, \  \ \text{ where } \
 c= \inf_{\vert x\vert \geq s_G,\, y \geq\mu_G} D_G(x-iy).
\]
The right-hand side defines an integrable function of $y$ over $[\eps, +\infty[$ : indeed, by Fubini's  Theorem, we find
\[
\int_\eps^{+\infty} \int_{\R^+}    | t G(t)| e^{-yt} dt dy = \int_{\R^+}  | G(t)| e^{-\eps t}dt.  
\]
The dominated convergence theorem allows us to conclude : we send $\mu$ and then $s$ to $+\infty$, and, like for Lemma \ref{TECHNIK}, we get that the second term vanishes. The same holds for the fourth term.  Finally, we obtain
\[
N_\epsilon= \frac1{2\pi i}\int_{\R}\frac{\widehat{\widetilde G}'(x-i\epsilon)}{D_G(x-i\epsilon)}dx,
\]
that is to say, $N_\epsilon$ is nothing but the winding number of the closed path defined by (continuity by)
\[
\left\{1-\int_{\R^+}G(t)e^{-i\omega t-\epsilon t}dt\right\}_{\omega\in \overline{\R}}.
\] 
Now, using that $\int_{\R^+}t|G(t)|dt<+\infty$ and $G\in L^1(\R^+)$, we obtain
\[
\lim_{\epsilon\to 0}\left\|\frac{\widehat{\widetilde G}'(\cdot-i\epsilon)}{1-\widehat{\widetilde G}(\cdot-i\epsilon)}-\frac{\widehat{\widetilde G}'(\cdot)}{1-\widehat{\widetilde G}(\cdot)}\right\|_{L^\infty(\R^+)}=0,
\]
from where we conclude $N_\epsilon=\text{Ind}_G(0)$ for all sufficiently small $\epsilon$. 
Therefore, under the assumption $D_G|_{\R}\neq 0$, the condition $D_G|_{\Pi^-}\neq 0$ is equivalent to requiring $\text{\rm Ind}_G(0)=0$.
 \hfill $\Box$

\section{Bootstrap argument, proof of Theorem \ref{MAINRES}}\label{S-PROOF}
Corollary \ref{FIRSTINEQ} indicates that, in order to get the conclusion of Theorem \ref{MAINRES}, it suffices to control the polynomial decay of the input term $F$ in equation \eqref{FULLEQ}. This control follows from a bootstrap argument that  involves appropriate Sobolev norms of the solution.  It is expressed in the next statement below. 

Given $n\in\N$, a solution $p$ of equation \eqref{GALILEANVAR} and $T>0$, consider the quantity $M_{n,T}(p)$
\[
M_{n,T}(p)=\max\left\{\sup_{t\in [0,T]}(1+t)^n|R(t)|,\sup_{t\in [0,T]}\frac{\|p(t)\|_{{\cal H}^n}}{1+t},\sup_{t\in [0,T]}\|p(t)\|_{{\cal H}^{n-2}}\right\}.
\]
\begin{Pro}
Assume that $g$ is of class $C^n(\R)$ for some $n\geq 4$ and satisfies the conditions \eqref{CONDG} and assume that the condition \eqref{STABCRIT} holds. There exists $M_{K}>0$, and for every $M\geq M_{K}$, there exists $\epsilon_{K,M}>0$ such that for every $\epsilon\in (0,\epsilon_{K,M})$ and every initial condition $p(0)= r(0)\cdot g$ with $\|p(0)\|_{{\cal H}^n}=1$, and for every $T>0$, whenever the inequality 
\[
M_{n,T}(p)\leq M,
\]
holds for the  corresponding solution of equation \eqref{GALILEANVAR}, we actually have $M_{n,T}(p)\leq \frac{M}{2}$. 
\label{BOOTSTRAP}
\end{Pro}
{\sl Proof of Theorem \ref{MAINRES}.} For a (rescaled) initial condition $\|r(0)\cdot g\|_{{\cal H}^n}=1$ as in Proposition \ref{BOOTSTRAP}, we have 
\[
M_{n,0}(p)\leq \max\{|R(0)|,1\}\leq \pi\sqrt{2},
\]
where the second inequality follows from the Cauchy-Schwarz inequality. Let $M_{\rm thr}>\max\{M_K,\pi\sqrt{2}\}$ 
and let $\epsilon_K=\epsilon_{K,M_{\rm thr}}$. We have $M_{n,T}(p)\leq M_{\rm thr}$ for $T>0$ sufficiently small, by continuity. Let $T_{\rm max}$ be defined as follows  
\[
T_{\rm max}=\sup\{T\in\R^+\ :\ M_{n,T}(p) \leq M_{\rm thr}\}.
\]
We claim that, for every $\epsilon\in (0,\epsilon_K)$, we have $T_{\rm max}=+\infty$ from where the theorem immediately follows.
In fact, if we had $T_{\rm max}<+\infty$, then Proposition \ref{BOOTSTRAP} would imply $M_{n,T_{\rm max}}(p) \leq \frac{M_{\rm thr}}{2}$. By continuity, there would exist $T>T_{\rm max}$ such that $M_{n,T}(p) \leq M_{\rm thr}$. But this  contradicts the definition of $T_{\rm max}$, hence we must have $T_{\rm max}=+\infty$. \hfill $\Box$

\subsection{Proof of Proposition \ref{BOOTSTRAP}}
We prove separately each of the three claims of the Proposition. Throughout the proof, the dependence on $K$ and $\epsilon$ is explicitly detailed so that Remark \ref{LARGEPERT} can be readily proved afterwards. 
The first step consists in propagating the estimate on $\sup_{t\in [0,T]}(1+t)^n|R(t)|$.
\begin{Lem}
Assume that $g$ is of class $C^n(\R)$ for some $n\geq 3$ and satisfies the conditions \eqref{CONDG} and assume that the condition \eqref{STABCRIT} holds. There exists $M_1>0$, and for every $M\geq M_1$, there exists $\epsilon_1>0$ so that for every $\epsilon\in (0,\epsilon_1)$ and every initial condition $p(0)$ with $\|p(0)\|_{{\cal H}^n}= 1$, and for every $T>0$, whenever the inequality 
\[
M_{n,T}(p)\leq M,
\]
holds for the subsequent solution of equation \eqref{GALILEANVAR}, we actually have $\sup\limits_{t\in [0,T]}(1+t)^n|R(t)|\leq \frac{M}{2}$.
\label{SECONDLEM}
\end{Lem}
{\sl Proof of the Lemma.}
Using that $p(s)$ is $n$-times differentiable at all times $s\in\R^+$, a reasoning similar to the one leading to the inequality \eqref{CONTFOUR} yields  
\[
\tau^{j}|\widehat{p}_k(s,\tau)|\leq \pi\sqrt{2}\|p(s)\|_{{\cal H}^j},\ \forall j\in\{0,\cdots,n\}, k\in\Z, s,\tau\in\R^+.
\]
Multiplying by ${n\choose j}$, summing for $j=0$ to $n$, we get 
\[
\sup_{t\in\R^+}(1+t)^n|\widehat{p}_k(s,t)|\leq 2^{n}\pi\sqrt{2}\|p(s)\|_{{\cal H}^n},\ \forall k\in\Z, s,\tau\in\R^+.
\]

Now, let $T>0$ be arbitrary and assume $\sup\limits_{t\in [0,T]}(1+t)^n|R(t)|\leq M$ for some $M>0$. Using Corollary \ref{FIRSTINEQ} together with the expression \eqref{FULLEQ} of $F$, the previous estimate and $\|p(0)\|_{{\cal H}^n}=1$, we successively have 
\begin{align}
\sup_{t\in [0,T]}(1+t)^n|R(t)|&\leq C_K\left(2^{n}\pi\sqrt{2}+\frac{\epsilon K}{2}M \sup_{t\in [0,T]}(1+t)^n\int_0^t\frac{|\widehat{p}_2(s,t+s)|}{(1+s)^n}ds\right)\nonumber\\
&\leq 2^{n}\pi\sqrt{2}C_K\left(1+\frac{\epsilon K}{2}M \sup_{t\in [0,T]}(1+t)^n\int_0^t\frac{\|p(s)\|_{{\cal H}^n}}{(1+t+s)^n (1+s)^n}ds\right)\nonumber\\
&\leq 2^{n}\pi\sqrt{2}C_K\left(1+\frac{\epsilon K}{2}M \sup_{t\in [0,T]}\int_0^t\frac{\|p(s)\|_{{\cal H}^n}}{(1+s)^n}ds\right)\nonumber\\
&\leq 2^{n}\pi\sqrt{2}C_K\left(1+\frac{\epsilon K}{2(n-2)}M^2\right)\label{DEFMEPS}
\end{align}
where the last inequality uses the assumption $\sup\limits_{s\in [0,t]}\frac{\|p(s)\|_{{\cal H}^n}}{1+s}\leq M$ for every $t\leq T$ and $\int_{\R^+}\frac{ds}{(1+s)^{n-1}}=\frac1{n-2}$. We conclude that the Lemma holds with 
\[
M_1=2^{n+2}\pi\sqrt{2}C_K\quad\text{and}\quad \epsilon_1=\frac{n-2}{2^{n+1}\pi\sqrt{2}C_KKM},
\]
(so that $2^{n}\pi\sqrt{2}C_K\leq \frac{M}4$ when $M_1\leq M$, and $2^{n}\pi\sqrt{2}C_K\frac{\epsilon K}{2(n-2)}M\leq \frac14$ when $\epsilon\leq \epsilon_1$). \hfill $\Box$

In order to propagate the bounds on the norms $\|p(t)\|_{{\cal H}^n}$ and $\|p(t)\|_{{\cal H}^{n-2}}$, we establish the following property. 
\begin{Lem}
Given $\ell\geq 1$, there exists a constant $C'_\ell>0$ such that for every $K\in\R^+$, we have 
\[
\frac{d\|p\|_{{\cal H}^{\ell}}}{dt}\leq C'_\ell K|R(t)|\left((1+t)^{\ell}(\|g\|_{{\cal H}^{\ell}}+\epsilon\|p\|_{{\cal H}^{0}})+\epsilon (1+t)\sum_{j=1}^{\ell}t^{\ell-j}\|p\|_{{\cal H}^{j}}\right), \forall t>0.
\]
\label{DIFFINEQ}
\end{Lem}
Strictly speaking, the inequality here applies to trajectories issued from smooth initial conditions, so that $t\mapsto \|p\|_{{\cal H}^{\ell}}$ is certainly differentiable. However, by a density argument, any inequality that follows suit from integration in time holds for trajectories in $C^n$ and this is what matters for the proofs of Lemma \ref{THIRDESTIM} and \ref{LASTLEMM} below.

\noindent
{\sl Proof of the Lemma.} 
Given the definition of $\|\cdot\|_{{\cal H}^\ell}$, all we need to control are the quantities $\frac{d\|\langle \omega\rangle\partial_\theta^{k_\theta}\partial_\omega^{k_\omega}p\|_{L^2(\T^1\times\R)}^2}{dt}$ for $k_\theta+k_\omega\leq \ell$. To that goal, using the scalar product associated with $\|\cdot\|_{L^2(\T^1\times\R)}$, we write
\[
\frac{d\|\langle \omega\rangle\partial_\theta^{k_\theta}\partial_\omega^{k_\omega}p\|_{L^2(\T^1\times\R)}^2}{dt}=2\int_{\T^1\times\R}\langle \omega\rangle^2\partial_t\partial_\theta^{k_\theta}\partial_\omega^{k_\omega}p\partial_\theta^{k_\theta}\partial_\omega^{k_\omega}pd\theta d\omega
\]
Now, by applying $\partial_\theta^{k_\theta}\partial_\omega^{k_\theta}$ to equation \eqref{GALILEANVAR}, we obtain that the equation for $\partial_t\partial_\theta^{k_\theta}\partial_\omega^{k_\omega}p=\partial_\theta^{k_\theta}\partial_\omega^{k_\omega}\partial_tp$ consists of three terms, namely 
\begin{itemize}
\item[$\bullet$] $\frac{1}{2\pi}\partial_\theta^{k_\theta}\partial_\omega^{k_\omega}\left(g(\omega)\partial_\theta W(\theta+t\omega,p)\right)=\frac{1}{2\pi}\partial_\omega^{k_\omega}\left(\partial_\theta^{k_\theta+1} W(\theta+t\omega,p)g(\omega)\right)$,
\item[$\bullet$] $\epsilon\partial_\theta^{k_\theta}\partial_\omega^{k_\omega}\left(\partial_\theta pW(\theta+t\omega,p)\right)$\item[$\bullet$] $\epsilon\partial_\theta^{k_\theta}\partial_\omega^{k_\omega}\left(p\partial_\theta W(\theta+t\omega,p)\right)$
\end{itemize}
which we analyze separately. 
To that goal we shall use the two basic properties. First, the partial derivative of the product $a\cdot b$ of two functions $a$ and $b$ of a real variable, say $x$, can be decomposed as follows 
\[
\partial_x^k(a\cdot b)=\sum_{j=0}^k{k\choose j}\partial_x^ja\cdot\partial_x^{k-j}b.
\]
Moreover, writing 
\[
W(\theta,p)=\frac{-iK}{2}\int_{\T^1\times\R}\left( e^{i(\theta'+t\omega'-\theta)}-c.c.\right)p(\theta',\omega')d\theta' d\omega',
\]
we easily compute for arbitrary integers $j_\theta,j_\omega$
\begin{align*}
\partial_\theta^{j_\theta}\partial_\omega^{j_\omega}W(\theta+t\omega,p)&=\frac{-iKt^{j_\omega}}{2}\int_{\T^1\times\R}\left((-i)^{j_\theta+j_\omega}e^{i(\theta'-\theta+t(\omega'-\omega))}-c.c.\right)p(\theta',\omega')d\theta' d\omega'\\
&=\frac{-iKt^{j_\omega}}{2}\left((-i)^{j_\theta+j_\omega}e^{-i(\theta+t\omega)}\overline{R(t)}-c.c.\right).
\end{align*}

$\bullet$ For the first term, we combine the two previous properties to obtain 
\begin{align*}
\partial_\omega^{k_\omega}\left(\partial_\theta^{k_\theta+1} W(\theta+t\omega,p)g(\omega)\right)&=\sum_{j_\omega=0}^{k_\omega}{k_\omega\choose j_\omega}\partial_\theta^{k_\theta+1}\partial_\omega^{j_\omega}W(\theta+t\omega,p)g^{(k_\omega-j_\omega)}(\omega)\\
&=\frac{-iK}{2}\sum_{j_\omega=0}^{k_\omega}{k_\omega\choose j_\omega}t^{j_\omega}\left((-i)^{k_\theta+j_\omega+1}e^{-i(\theta+t\omega)}\overline{R(t)}-c.c.\right)g^{(k_\omega-j_\omega)}(\omega).
\end{align*}
Multiplying by $2\langle \omega\rangle^2\partial_\theta^{k_\theta}\partial_\omega^{k_\omega}p$ and integrating over $\T^1\times\R$, this expression gives the following contribution to $\frac{d\|\langle \omega\rangle\partial_\theta^{k_\theta}\partial_\omega^{k_\omega}p\|_{L^2(\T^1\times\R)}^2}{dt}$
\begin{align*}
\frac{1}{\pi}\int_{\T^1\times\R}\langle \omega\rangle^2\partial_\omega^{k_\omega}\left(\partial_\theta^{k_\theta+1} W(\theta+t\omega,p)g(\omega)\right)\partial_\theta^{k_\theta}\partial_\omega^{k_\omega}pd\theta d\omega=\\
\frac{-iK}{2\pi}\int_{\T^1\times\R}\left(\langle \omega\rangle^2\sum_{j_\omega=0}^{k_\omega}{k_\omega\choose j_\omega}t^{j_\omega}\left((-i)^{k_\theta+j_\omega+1}e^{-i(\theta+t\omega)}\overline{R(t)}-c.c.\right)g^{(k_\omega-j_\omega)}(\omega)\right)\partial_\theta^{k_\theta}\partial_\omega^{k_\omega}pd\theta d\omega
\end{align*}
Using $\left|-i\left((-i)^{k_\theta+j_\omega+2}e^{-i(\theta+t\omega)}\overline{R(t)}-c.c.\right)\right|\leq 2|R(t)|$ and the Cauchy-Schwarz inequality, this expression turns out to be bounded above by 
\[
\frac{K}{\pi}|R(t)|\|\langle \omega\rangle\partial_\theta^{k_\theta}\partial_\omega^{k_\omega}p\|_{L^2(\T^1\times\R)}\sum_{j_\omega=0}^{k_\omega}{k_\omega\choose j_\omega}t^{j_\omega}\|\langle\omega\rangle g^{(k_\omega-j_\omega)}\|_{L^2(\T^1\times\R)}.
\]
By summing over $k_\theta$ and $k_\omega$, and using the inequality (straightforward consequence of Cauchy-Schwarz)
\[
\sum_{k=1}^m|a_k|\leq 2^{m/2}\left(\sum_{k=1}^m|a_k|^2\right)^{1/2},
\]
for every $m\in\N$ and every $\{a_k\}_{k=1}^m\in\C^m$, we conclude that the first term of the differential inequality for $\frac{d\|p(t)\|_{{\cal H}^\ell}^2}{dt}$ is of the form
\[
CK(1+t)^\ell|R(t)|\|g\|_{{\cal H}^{\ell}}\|p\|_{{\cal H}^{\ell}},
\]
for some constant $C>0$.

$\bullet$ For the second term, we first use the basic properties above to get 
\begin{align*}
\partial_\omega^{k_\omega}(\partial_\theta pW(\theta+t\omega,p))&=\sum_{j_\omega=0}^{k_\omega}{k_\omega\choose j_\omega}\partial_\theta\partial_\omega^{j_\omega}p\partial_\omega^{k_\omega-j_\omega}W(\theta+t\omega,p)\\
&=\frac{-iK}{2}\sum_{j_\omega=0}^{k_\omega}c_{j_\omega,k_\omega}t^{k_\omega-j_\omega}\left((-i)^{k_\omega-j_\omega}e^{-i(\theta+t\omega)}\overline{R(t)}-c.c.\right)\partial_\theta\partial_\omega^{j_\omega}p
\end{align*}
and then
\begin{multline*}
\partial_\theta^{k_\theta}\partial_\omega^{k_\omega}(\partial_\theta pW(\theta+t\omega,p))=
\\
\frac{-iK}{2}\sum_{j_\theta=0}^{k_\theta}\sum_{j_\omega=0}^{k_\omega}c_{j_\theta,k_\theta}c_{j_\omega,k_\omega}t^{k_\omega-j_\omega}\left((-i)^{k_\theta-j_\theta+k_\omega-j_\omega}e^{-i(\theta+t\omega)}\overline{R(t)}-c.c.\right)\partial_\theta^{j_\theta+1}\partial_\omega^{j_\omega}p.
\end{multline*}
We consider separately the term $j_\omega=k_\omega$. For $j_\theta=k_\theta$, by multiplying by $\partial_\theta^{k_\theta}\partial_\omega^{k_\omega}p$ and integrating over $\T^1$, we obtain
\begin{align*}
\int_{\T^1}\left(e^{-i(\theta+t\omega)}\overline{R(t)}-c.c.\right)\partial_\theta^{k_\theta+1}\partial_\omega^{k_\omega}p\partial_\theta^{k_\theta}\partial_\omega^{k_\omega}pd\theta =\frac{1}2\int_{\T^1}\left(e^{-i(\theta+t\omega)}\overline{R(t)}-c.c.\right)\partial_\theta(\partial_\theta^{k_\theta}\partial_\omega^{k_\omega}p)^2d\theta\\
=-\frac{1}2\int_{\T^1}\partial_\theta\left(e^{-i(\theta+t\omega)}\overline{R(t)}-c.c.\right)(\partial_\theta^{k_\theta}\partial_\omega^{k_\omega}p)^2 d\theta\\
=-\frac{1}2\int_{\T^1}\left(-ie^{-i(\theta+t\omega)}\overline{R(t)}-c.c.\right)(\partial_\theta^{k_\theta}\partial_\omega^{k_\omega}p)^2d\theta
\end{align*}
where the second equality follows from integration by parts and periodicity. Multiplying by $\langle \omega\rangle$ and integrating over $\R$, it follows that
\[
-i\int_{\T^1\times\R}\langle \omega\rangle^2\left(e^{-i(\theta+t\omega)}\overline{R(t)}-c.c.\right)\partial_\theta^{k_\theta+1}\partial_\omega^{k_\omega}p\partial_\theta^{k_\theta}\partial_\omega^{k_\omega}pd\theta d\omega\leq |R(t)|\|\langle \omega\rangle\partial_\theta^{k_\theta}\partial_\omega^{k_\omega}p\|_{L^2(\T^1\times\R)}^2.
\]
For the remaining terms $j_\theta<k_\theta$ and $j_\omega=k_\omega$, using Cauchy-Schwarz inequality again, we obtain
\begin{multline*}
-i\int_{\T^1\times\R}\langle \omega\rangle^2\left((-i)^{k_\theta-j_\theta}e^{-i(\theta+t\omega)}\overline{R(t)}-c.c.\right)\partial_\theta^{j_\theta+1}\partial_\omega^{k_\omega}p\partial_\theta^{k_\theta}\partial_\omega^{k_\omega}pd\theta d\omega
\leq 
\\
|R(t)|\|\langle \omega\rangle\partial_\theta^{j_\theta+1}\partial_\omega^{k_\omega}p\|_{L^2(\T^1\times\R)}
\|\langle \omega\rangle\partial_\theta^{k_\theta}\partial_\omega^{k_\omega}p\|_{L^2(\T^1\times\R)},
\end{multline*}
and $j_\theta+1\leq k_\theta$. Altogether, by summing over $k_\theta,k_\omega$ and $j_\theta$, we obtain that the terms $j_\omega=k_\omega$ give a total contribution to the differential inequality for $\partial_t\|p(t)\|_{{\cal H}^\ell}$ of the form
\[
CK|R(t)|\|p\|_{{\cal H}^{\ell}}^2,
\]
where $C>0$ is again a generic constant. 

For the terms $j_\omega<k_\omega$, ignoring the constants and the factors involving $\overline{R(t)}$, we consider the following change of index
\[
\sum_{j_\omega=0}^{k_\omega-1}t^{k_\omega-j_\omega}\partial_\theta^{j_\theta+1}\partial_\omega^{j_\omega}p=\sum_{j_\omega=1}^{k_\omega}t^{j_\omega}\partial_\theta^{j_\theta+1}\partial_\omega^{k_\omega-j_\omega}p.
\]
In this expression, the derivative indices satisfy the inequality $j_\theta+1+k_\omega-j_\omega\leq \ell-j_\omega+1$. Hence, we can repeat the same procedure as before. After summation over $j_\theta,k_\theta$ and $k_\omega$, we get that the terms $j_\omega<k_\omega$ yield a contribution to the differential inequality for $\frac{d\|p(t)\|_{{\cal H}^\ell}^2}{dt}$ of the form
\[
CK|R(t)|\sum_{j=1}^{\ell}t^{j}\|p\|_{{\cal H}^{\ell-j+1}}\|p\|_{{\cal H}^{\ell}}.
\]

$\bullet$ The computation for the third term is similar to that of the second term and results in a total contribution of the form 
\[
CK|R(t)|\sum_{j=0}^{\ell}t^{j}\|p\|_{{\cal H}^{\ell-j}}\|p\|_{{\cal H}^{\ell}}.
\]
Adding all contributions together and using that $\frac{d\|p(t)\|_{{\cal H}^\ell}^2}{dt}=2\|p\|_{{\cal H}^{\ell}}\partial_l\|p\|_{{\cal H}^{\ell}}$, the conclusion of the Lemma easily follows. \hfill $\Box$ 

We can now pass to the propagation of the estimates on $\frac{\|p(t)\|_{{\cal H}^n}}{1+t}$ and $\|p(t)\|_{{\cal H}^{n-2}}$.
\begin{Lem}
Assume that $g$ is of class $C^n(\R)$ for some $n\geq 4$ and satisfies the conditions \eqref{CONDG} and assume that the condition \eqref{STABCRIT} holds. There exists $M_2>0$, and for every $M\geq M_2$, there exists $\epsilon_2>0$ such that for every $\epsilon\in (0,\epsilon_{2})$ and every initial condition $p(0)$ with $\|p(0)\|_{{\cal H}^n}= 1$, and for every $T>0$, whenever the inequality 
\[
M_{n,T}(p)\leq M,
\]
holds for the subsequent solution of equation \eqref{GALILEANVAR}, we actually have $\sup\limits_{t\in [0,T]}\frac{\|p(t)\|_{{\cal H}^n}}{1+t}\leq \frac{M}{2}$.
\label{THIRDESTIM}
\end{Lem}
{\sl Proof of the Lemma.} Using the estimate \eqref{DEFMEPS} from the proof of Lemma \ref{SECONDLEM} together with Lemma \ref{DIFFINEQ} for $\ell=n$ (and using the notation $C_2=2^{n}\pi\sqrt{2}C'_{n}C_K$),
we obtain
\begin{align*}
\frac{d\|p\|_{{\cal H}^n}}{dt}&\leq C_2K\left(1+\frac{\epsilon K}{2(n-2)}M^2\right)\left(\|g\|_{{\cal H}^{n}}+\epsilon\|p\|_{{\cal H}^{0}}+\epsilon \sum_{j=1}^{n}\frac{\|p\|_{{\cal H}^{j}}}{(1+t)^{j-1}}\right)\\
&\leq C_2K\left(1+\frac{\epsilon K}{2(n-2)}M^2\right)\left(\|g\|_{{\cal H}^{n}}+\epsilon(n-1)M+2\epsilon \frac{\|p\|_{{\cal H}^{n}}}{(1+t)^{n-2}}\right)
\end{align*}
where the second inequality relies both on $\|p\|_{{\cal H}^{j}}\leq \|p\|_{{\cal H}^{n-2}}$ for $j\in\{0,\cdots,n-2\}$ and on $\|p\|_{{\cal H}^{n-1}}\leq \|p\|_{{\cal H}^{n}}$. Using $\|p(0)\|_{{\cal H}^{n}}= 1$, integration and the Gronwall inequality then successively yield
\begin{align*}
\frac{\|p(t)\|_{{\cal H}^n}}{1+t}&\leq \max\left\{1,C_2K\left(1+\frac{\epsilon K}{2(n-2)}M^2\right)\left(\|g\|_{{\cal H}^{n}}+\epsilon(n-1)M\right)\right\}\\
&+\epsilon C_2K\left(1+\frac{\epsilon K}{2(n-2)}M^2\right)\int_0^t\frac{\|p(s)\|_{{\cal H}^{n}}}{(1+t)(1+s)^{n-2}}ds\\
&\leq \max\left\{1,C_2K\left(1+\frac{\epsilon K}{2(n-2)}M^2\right)\left(\|g\|_{{\cal H}^{n}}+\epsilon(n-1)M\right)\right\}e^{\frac{\epsilon C_2K\left(1+\frac{\epsilon K}{2(n-2)}M^2\right)}{n-3}},
\end{align*}
for all $t\in [0,T]$. By evaluating this quantity for $\epsilon=0$ and using monotonicity with respect to $\epsilon$, we conclude that the Lemma holds with 
\[
M_2=8\max\{C_2K\|g\|_{{\cal H}^{n}},1\},
\]
and $\epsilon_2$ being the largest $\epsilon>0$ such that we simultaneously have 
\[
e^{\frac{\epsilon C_2K\left(1+\frac{\epsilon K}{2(n-2)}M^2\right)}{n-4}}\leq 2,
\]
and
\[
\epsilon C_2K\left(n-1+\frac{K}{2(n-2)}M\|g\|_{{\cal H}^{n}}+\frac{\epsilon K(n-1)}{2(n-2)}M^2\right)\leq \frac18.
\]
\hfill $\Box$

Finally, we proceed similarly to propagate the estimate on $\|p(t)\|_{{\cal H}^{n-2}}$.
\begin{Lem}
Assume that $g$ is of class $C^n(\R)$ for some $n\geq 2$ and satisfies the conditions \eqref{CONDG} and assume that the condition \eqref{STABCRIT} holds. There exists $M_3>0$, and for every $M\geq M_3$, there exists $\epsilon_3>0$ such that for every $\epsilon\in (0,\epsilon_{3})$ and every initial condition $p(0)$ with  and $\|p(0)\|_{{\cal H}^n}= 1$, and for every $T>0$, whenever the inequality 
\[
M_{n,T}(p)\leq M,
\]
holds for the subsequent solution of equation \eqref{GALILEANVAR}, we actually have $\sup\limits_{t\in [0,T]}\|p(t)\|_{{\cal H}^{n-2}}\leq \frac{M}{2}$.
\label{LASTLEMM}
\end{Lem}
{\sl Proof of the Lemma.} Proceeding similarly as in the previous proof, we obtain
\[
\frac{d\|p\|_{{\cal H}^{n-2}}}{dt}\leq C_3K\left(1+\frac{\epsilon K}{2(n-2)}M^2\right)\left(\|g\|_{{\cal H}^{n-2}}+\epsilon(n-1)M\right)\frac1{(1+t)^2},
\]
where $C_3=2^{n}\pi\sqrt{2}C'_{n-2}C_K$. Using that $\int_{\R^+}\frac{dt}{(1+t)^2}=1$, we then get after integration (using also $\|r\|_{{\cal H}^{n-2}}\leq \|r\|_{{\cal H}^{n}}$)
\[
\|p(t)\|_{{\cal H}^{n-2}}\leq 1+C_3K\left(1+\frac{\epsilon K}{2(n-2)}M^2\right)\left(\|g\|_{{\cal H}^{n-2}}+\epsilon(n-1)M\right)
\]
from where the lemma follows  with 
\[
M_3=4(1+C_3K\|g\|_{{\cal H}^{n-2}}),
\]
and $\epsilon_3$ defined as the largest $\epsilon>0$ such that 
\[
\epsilon C_3K\left(n-1+\frac{\epsilon K}{2(n-2)}M\|g\|_{{\cal H}^{n-2}}+\epsilon \frac{\epsilon K(n-1)}{2(n-2)}M^2\right)= \frac14.
\]
\hfill $\Box$

\noindent
The proposition finally holds with $M_{K}=\max\left\{M_1,M_2,M_3\right\}$ and $\epsilon_{K,M}=\min\left\{\epsilon_1,\epsilon_2,\epsilon_3\right\}$.

\subsection{Proof of Corollary \ref{WEAKCONV}}\label{S-WEAKCONV}
We aim at showing that the solution of equation \eqref{GALILEANVAR}, which by Proposition \ref{BOOTSTRAP} satisfies the bound $\sup\limits_{t\in\R^+}\|p(t)\|_{{\cal H}^{n-2}}<+\infty$, converges in ${\cal H}^{n-2}$ to the function $p_\infty$ defined by
\[
p(0,\theta,\omega)+\int_{\R^+}\left(\epsilon\partial_\theta p(s,\theta,\omega)W(\theta+s\omega,p(s))+\left(\frac{g(\omega)}{2\pi}+\epsilon p(s,\theta,\omega)\right)\partial_\theta W(\theta+s\omega,p(s))\right)ds,
\]
for all $(\theta,\omega)\in\T^1\times\R$. (NB: The proof simultaneously shows that $p_\infty$ is well-defined in ${\cal H}^{n-2}$.)

To that goal, it suffices to control the quantity $\|p(t)-p_\infty\|_{{\cal H}^{n-2}}^2$, hence to control each of the integrals
\[
I_1(t)=\int_t^{+\infty}\|\partial_\theta p(s)W(\theta+s\omega,p(s))\|_{{\cal H}^{n-2}}^2ds,\quad
\int_t^{+\infty}\|g(\omega)\partial_\theta W(\theta+s\omega,p(s))\|_{{\cal H}^{n-2}}^2ds,
\]
and
\[
\int_t^{+\infty}\| p(s)\partial_\theta W(\theta+s\omega,p(s))\|_{{\cal H}^{n-2}}^2ds.
\]
Using the estimates obtained for each term in the proof of Lemma \ref{DIFFINEQ} above, one obtains the following inequalities
\begin{align*}
I_1(t)&\leq CK\int_t^{+\infty}|R(s)|\left(\|p(s)\|_{{\cal H}^{n-2}}+\sum_{j=1}^{n-2}t^j\|p(s)\|_{{\cal H}^{n-1-j}}\right)\|p(s)\|_{{\cal H}^{n-2}}ds\\
&\leq  CK M_{\rm thr}^3(n-2)\int_t^{+\infty}\frac{ds}{(1+s)^2},
\end{align*}
(where the second inequality relies on Proposition \ref{BOOTSTRAP}), and similar inequalities hold for the two other integrals. The asymptotic behavior $\lim\limits_{t\to+\infty}\|p(t)-p_\infty\|_{{\cal H}^{n-2}}^2=0$ then immediately follows and Corollary \ref{WEAKCONV} is proved.

\subsection{Proof of Proposition \ref{LARGEPERT}}\label{S-LARGEPERT}
To prove the Proposition, beside observing that the stability critierion \eqref{STABCRIT} can always be satisfied by choosing $K$ sufficiently small, it suffices to establish a statement analogous to Proposition \ref{BOOTSTRAP} in which the roles of $K$ and $\epsilon$ are exchanged. Accordingly, one has to verify that $K$ and $\epsilon$ can be exchanged in the statements of Lemmas \ref{SECONDLEM}, \ref{THIRDESTIM} and \ref{LASTLEMM}. 

For Lemma \ref{SECONDLEM}, this is immediate from expression \eqref{DEFMEPS}. For Lemma \ref{THIRDESTIM}, we first notice that the constant $C_K$ in Corollary \ref{FIRSTINEQ} obviously depends continuously on $K$; thus so does the constant $C_2=C_2(K)$ in the proof of Lemma \ref{THIRDESTIM}. Clearly, the statement of Lemma \ref{THIRDESTIM} holds for every $\epsilon>0$, $M\geq 4$ and $K\leq K_{\epsilon,M}$ where $K_{\epsilon,M}>0$ is sufficiently small so that we simultaneously have
\[
e^{\frac{\epsilon C_2(K_{\epsilon,M})K_{\epsilon,M}\left(1+\epsilon c_{n-2} K_{\epsilon,M}M^2\right)}{n-4}}\leq 2,
\]
and
\[
\epsilon C_2(K_{\epsilon,M})K_{\epsilon,M}\left(n-1+c_{n-2}K_{\epsilon,M}M\|g\|_{{\cal H}^{n}}+\epsilon c_{n-2}(n-1)K_{\epsilon,M}M^2\right)\leq \frac18.
\]
The reasoning is similar for Lemma \ref{LASTLEMM} and the proof of Proposition \ref{LARGEPERT} is complete. 

\noindent
{\bf Acknowledgements} 

\noindent
Work supported by CNRS PEPS "Physique Th\'eorique et ses Interfaces".

\end{document}